\documentclass[12pt,leqno]{article}
\usepackage{amsfonts}
\usepackage{amsmath, amsthm, amsfonts, amssymb, color}
\usepackage{mathrsfs}
\oddsidemargin=0mm \topmargin=-12mm \numberwithin{equation}{section}

\newcommand{\be}{\begin{eqnarray}}
\newcommand{\ee}{\end{eqnarray}}
\newcommand{\ce}{\begin{eqnarray*}}
\newcommand{\de}{\end{eqnarray*}}
\newtheorem{theorem}{Theorem}[section]
\newtheorem{lemma}[theorem]{Lemma}
\newtheorem{remark}[theorem]{Remark}
\newtheorem{definition}[theorem]{Definition}
\newtheorem{proposition}[theorem]{Proposition}
\newtheorem{Examples}[theorem]{Example}
\newtheorem{corollary}[theorem]{Corollary}

\def\v{{\mathrm{v}}}
\def\eps{\varepsilon}

\def\p{\partial}
\def\d{\delta}

\def\[{{\Big[}}
\def\]{{\Big]}}
\def\<{{\langle}}
\def\>{{\rangle}}
\def\({{\Big(}}
\def\){{\Big)}}

\def\bx{{\mathbf{x}}}

\def\dif{{\mathord{{\rm d}}}}

\def\no{\nonumber}
\def\={&\!\!=\!\!&}
\def\bt{\begin{theorem}}
\def\et{\end{theorem}}
\def\bl{\begin{lemma}}
\def\el{\end{lemma}}
\def\br{\begin{remark}}
\def\er{\end{remark}}

\def\bd{\begin{definition}}
\def\ed{\end{definition}}
\def\bp{\begin{proposition}}
\def\ep{\end{proposition}}
\def\bc{\begin{corollary}}
\def\ec{\end{corollary}}
\def\bx{\begin{Examples}}
\def\ex{\end{Examples}}

\def\cB{{\mathcal B}}

\def\cJ{{\mathcal J}}

\def\cL{{\mathcal L}}

\def\beq{\beg{equation}}
\def\ff{\frac}

\def\mE{{\mathbb E}}

\def\mH{{\mathbb H}}

\def\mR{{\mathbb R}}
\def\mS{{\mathbb S}}

\def\mU{{\mathbb U}}

\def\mW{{\mathbb W}}

\def\sB{{\mathscr B}}
\def\sC{{\mathscr C}}
\def\sD{{\mathscr D}}

\def\sF{{\mathscr F}}

\def\E{\mathbb E}
\def\beg{\begin}

\def\geq{\geqslant}
\def\leq{\leqslant}
\def\ge{\geqslant}
\def\le{\leqslant}

\def\R{\mathbb R} \def\B{\mathcal B}

\def\nn{\nabla}\def\si{\sigma}\def\DD{\Delta}\def\dd{\delta}\def\bb{\beta}\def\aa{\alpha}\def\gg{\gamma}
\def\vv{\varepsilon}\def\e{{\rm e}}

\def\ss{\sqrt}\def\GG{\Gamma}\def\d{\dif}

\allowdisplaybreaks

\begin{document}

\title{{\bf Gradient Estimates for SDEs Driven by Multiplicative L\'evy Noise}
\footnote{Supported in part by NNSFC(11131003, 11431014, 11271294, 11325105), Macau 
S.A.R FDCT 049/2014/A1, Program for New Century Excellent Talents in University, the 985 project and the Laboratory of Mathematical and  Complex Systems.}}

\date{}
\author{
{\bf Feng-Yu Wang$^{a),b)}$\footnote{Corresponding author: wangfy@bnu.edu.cn; F.Y.Wang@swansea.ac.uk} Lihu Xu$^{c)}$, Xicheng  Zhang$^{d), e)}$ }\\
\footnotesize{$^{a)}$School of Mathematical Sciences, Beijing Normal
University, Beijing 100875, China}\\
\footnotesize{$^{b)}$Department of Mathematics,
Swansea University, Singleton Park, SA2 8PP, UK}\\
 \footnotesize{$^{c)}$Faculty of Science and Technology,
University of Macau, Av. Padre Tom¨¢s Pereira, }\\
 \footnotesize{ Taipa Macau, China}\\
\footnotesize{$^{d)}$  School of Mathematics and Statistics,
Wuhan University, Wuhan 430072, China}\\
\footnotesize{$^{e)}$Computational Science Hubei Key Laboratory, Wuhan University, Wuhan, 430072, China} }

 \maketitle
\begin{abstract} Gradient estimates are derived, for the first time,  for the semigroup associated to a class of stochastic differential equations driven by multiplicative L\'evy noise.  In particular, the  estimates are sharp for $\aa$-stable type noises.  To derive these estimates,
a new derivative formula of Bismut-Elworthy-Li's type
is established for the semigroup by using the Malliavin calculus and a finite-jump approximation argument.

\

 \end{abstract}

\section{Introduction}

Consider the following  stochastic differential equation (abbreviated as SDE) in $\mR^d$:
\begin{align}
\label{E0} \dif X_t =  b_t(X_t) \dif t +\si_t(X_{t-})\dif L_t,\ \ X_0=x,
\end{align}
where $L_t$ is a $d$-dimensional L\'evy process,   and
$$b: [0,\infty)\times \mR^d\to \mR^d,\ \ \si: [0,\infty)\times \mR^d\to \mR^d\otimes\mR^d$$ are   continuous  such that
($\nn_\v$ is the $\v$-directional derivative)
\begin{align}
\label{A1} \max\big\{\<\nn_\v b_s, \v\>,  \|  \nn_\v \si_s \|_{\mathrm{H.S.}}^2\big\}\leq K_s|\v|^2,\ \  s\ge 0,\ \v\in\mR^d
\end{align}
holds for some positive function
$K\in C([0,\infty);\mR_+).$  Under \eqref{A1}, it is well-known that
 (\ref{E0}) has a unique solution  $X_t(x)$ starting from $x$.    We aim to investigate the gradient estimate of $P_t:$
$$P_t f(x):= \mE f(X_t(x)),\ \ t\ge 0, f\in \cB_b(\mR^d),$$
where $\cB_b(\mR^d)$ denotes the space of all bounded Borel measurable functions on $\mR^d$.

When $b_t(x)=Ax$ for a matrix $A$ and $\si_t(x)=\text{Id}$, the gradient estimate of type
$$|\nn P_tf|\leq \|f\|_\infty \varphi(t),\ \ t>0,\ \ f\in \cB_b(\mR^d)$$ has been derived in \cite{W11}
by using lower bound conditions of the L\'evy measure,
and in \cite{SW} for $A=0$
by using
asymptotic behaviours of the symbol of $L_t$; see also \cite{W12} for a derivative formula by using coupling through the Mecke formula.

 Recently, a time-change argument was introduced in \cite{Zh1} to establish
 Bismut-Elworthy-Li's derivative formula of $P_t$
 for the case that
 $L_t$ is an $\aa$-stable process,
 $\si_t$ is invertible and independent of the space variable (i.e. the noise is additive), and $\nn b_t$ is uniformly bounded.  In particular, this derivative formula
implies   that for any $p>1$ there exists a constant $C(p)>0$ such that (see   \cite[Theorem 1.1]{Zh1})
\begin{align}\label{Z} |\nn P_t f|\leq \ff{C(p)}{1\land t^{\ff 1 \aa}} (P_t|f|^p)^{1/p},\ \ f\in\cB_b(\mR^d), t>0.\end{align} Using this time-change argument and the coupling method, Harnack inequalities are then established in \cite{WW}.

In this paper, we intend to extend the gradient estimate (\ref{Z}) for more general symmetric L\'evy process $L_t$   and space-dependent $\si_t$  (i.e. the noise is multiplicative).
Notice that if $b$ and $\sigma$ are independent of $t$, then the generator of $P_t$ is given by
$$
\cL f(x)=b(x)\cdot\nabla f(x)+\mathrm{P.V.}\int_{\mR^d}[f(x+\sigma(x)y)-f(x)]\nu(\dif y),
$$
where $\nu(\dif y)$ is the L\'evy measure of $L_t$, and P.V. stands for the Cauchy principal value.

From now on, we let $L_t$ be a symmetric and rotationally invariant L\'evy process, which can be formulated as
subordination of the Brownian motion.
More precisely,
let $W_t$ be a $d$-dimensional Brownian motion and $S_t$ an independent
subordinator
associated with a Bernstein function $B$ with $B(0)=0$; i.e. $S_t$ is an increasing process with stationary independent increments such that
$$
\mE \e^{-uS_t} = \e^{-t B(u)},\ \ u\geq 0, t\geq 0.
$$
Then $L_t:= W_{S_t}$ is a L\'evy process with symbol $\Psi(\xi):= B(|\xi|^2)$ (see e.g. \cite{J}).
In particular, if $B(u)=u^{\ff\aa 2}$ for some constant $\aa\in (0,2),$
then $L_t$ is an $\aa$-stable process.
By using the above relalization, equation (\ref{E0}) is now reduced to
\begin{align}\label{E1}
\dif X_t = b_t(X_t) \dif t +\si_t(X_{t-})\dif  W_{S_t}.
\end{align}
We assume

\

\begin{enumerate} \item[{\bf (A)}]  $\nabla b_t(x), \nabla\sigma_t(x)$ are locally (uniformly in $t$) Lipschitz continuous in $x$,
such that for any $t>0$,
\begin{equation}\begin{split} \label{Ass0}
&\|\nabla b\|_t:=\sup_{s\in [0,t], x\in\mR^d}\|\nn b_s(x)\|<+\infty,\\
&\|\nabla \sigma\|_t:= \sup_{s\in [0,t], x\in\mR^d}\|\nn\si_s(x)\|<+\infty,
\end{split}\end{equation}
 where $\|\cdot\|$ stands for the operator norm, and  for some positive increasing function $c_\cdot$ and some constant $m\geq 0$
\begin{align}
\sup_{s\in [0,t]} \|\sigma^{-1}_s(x)\|\leq c_t(1+|x|^m),\ \ t\geq 0, x\in\mR^d.\label{Ass1}
\end{align}
\end{enumerate}

\

Obviously, (\ref{Ass0}) implies (\ref{A1}) and hence, the existence and uniqueness of the solution.
Before stating our main result, let us briefly recall the main argument introduced in \cite{Zh1} for the study of the additive noise case, and explain why this argument is no longer valid in the multiplicative case.

In the additive noise case where $\si_t(x)=\si_t$ is independent of $x$,    for a fixed path $\ell_\cdot$ of $S_\cdot$,  we may reformulate the equation (\ref{E1}) as
$$\d X_t^\ell =b_t(X_t^\ell)\d t + \si_t\d  W_{\ell_t}.$$ To establish a derivative formula for $P_t^\ell$ with  $P_t^\ell f:= \mE f(X_t^\ell)$, consider
the following regularization of $\ell_\cdot$ for $\vv>0$:
$$\ell_t^\vv:= \ff 1\vv \int_t^{t+\vv} \ell_s\d s +\vv t,\ \ t\ge 0,$$ and the associated stochastic differential equation
$$ \d X_t^\vv = b_t(X_t^\vv)\d t +\si_t \d W_{\ell_t^\vv},\ \ X_0^\vv=X_0^\ell.$$
Then for any $t>0$, as $\vv\downarrow 0$ we have $\ell_t^\vv\downarrow \ell_t$ and $X_t^\vv\to X_t^\ell, \nn X_t^\vv\to\nn X_t^\ell$ in $L^p(\mathbb P)$ for any $p>1.$ Since $\ell_t^\vv$ is absolutely continuous, $X_t^\vv$ is indeed a diffusion process, so that an existing
Bismut-Elworthy-Li's type
derivative formula applies. Therefore,
letting $\vv\downarrow 0$ we derive a derivative formula for $P_t^\ell.$

Now, coming back to the multiplicative noise case, we consider
$$\d X_t^\ell =b_t(X_t^\ell)\d t + \si_t(X_{t-}^\ell)\d  W_{\ell_t},$$ and the corresponding approximation equation
$$ \d X_t^\vv = b_t(X_t^\vv)\d t +\si_t(X_{t}^\vv) \d W_{\ell_t^\vv},\ \ X_0^\vv=X_0^\ell.$$
Since we only have weak convergence $\d\ell_t^\vv \to \d\ell_t$ as Lebesgue-Stieltjes measures,
and the function $s\mapsto \si_s(X_{s-}^\ell)$ is  discontinuous,
the assertion that $X_t^\vv\to X_t^\ell$ is no longer true!  To see this, let us look at a  simple  example with $b_t(x)=0, \si_t(x)=\ss{1+|x|^2}$ and $d=1$. In this case, for   $X_0^\ell=0$ we have
$$\d |X_t^\vv|^2 = (1+|X_t^\vv|^2)\d\ell_t^\vv + 2 X_t^\vv \ss{1+|X_t^\vv|^2}\,\d W_{\ell_t^\vv},$$ so that
$\E(1+|X_t^\vv|^2)= \e^{\ell_t^\vv}.$ Then, if $X_t^\vv\to X^\ell_t$ as $\vv\downarrow 0$,  
by the dominated convergence theorem, we obtain
\beg{equation}\label{QB}\E(1+|X_t^\ell|^2)=\lim_{\vv\downarrow 0} \E(1+|X_t^\vv|^2)= \e^{\ell_t},\ \ \  t\ge 0.\end{equation}
On the other hand,
notice that $X_t^\ell= \int_0^t \ss{1+|X_{s-}^\ell|^2}\d W_{\ell_t}$ implies
$$\E(1+|X_t^\ell|^2)= 1+\int_0^t \E(1+|X_{s-}^\ell|^2)\d\ell_s,\ \ t\ge 0.$$ Combining this with \eqref{QB} we arrive at
$$
\e^{\ell_t}=1+\int^t_0 \e^{\ell_{s-}}\d\ell_s,\ \ t\ge 0,
$$ 
which is however not true for non-trivial $\ell_t$, since  $\DD \ell_{t_0}:=\ell_{t_0}-\ell_{t_0-}>0$  implies
$$
\e^{\ell_{t_0}}-\e^{\ell_{t_0-}}=\e^{\ell_{t_0-}}(\e^{\DD\ell_{t_0}}-1)> \e^{\ell_{t_0-}} \DD \ell_{t_0}.
$$

Due to the invalidity of regularization approximations,    we  have to  introduce a different argument.
Fortunately, by using a finite-jump approximation (i.e. approximating $\ell$ by those of finite many jumps in finite intervals) and the Malliavin calculus, we are able to establish a nice derivative formula for $P_t^\ell$ (see Theorem 3.2 below), which
in turn implies
the following main result of the paper.

\bt
Assume {\bf (A)} and let $P_t$ be the semigroup associate to equation $(\ref{E1})$.
\begin{enumerate}\item[{\rm (i)}]
  For fixed $R>0$, let $\tau:=\inf\{t: S_t\geq R\}$. If $\mE S_{t}^{-1/2}<\infty$, then for any $f\in\mathcal B_b(\mR^d)$, $\v\in\mR^d$ with $|\v|=1$, we have
\begin{align}
\nabla_{\mathrm{v}} P_tf=\mE\Bigg[& f(X_t)\frac{1}{S_{t\wedge\tau}}\Bigg(\int^{t\wedge\tau}_0
\<\sigma^{-1}_s(X_{s-})\nabla_\v X_{s-},\dif W_{S_s}\>\no\\
&-\int^{t\wedge\tau}_0\mathrm{Tr}\big(\sigma^{-1}_s\nabla_{\nabla_\v X_{s-}}\sigma_s\big) (X_{s-}) \dif S_s\label{EE11}\\
& +\int^{t\wedge\tau}_0\!\!\!\int_{\mR^d}\<\sigma^{-1}_s(X_{s-})\nabla_{\nabla_\v X_{s-}}\sigma_s(X_{s-})y,y\>N(\dif s,\dif y)\Bigg)\Bigg],
\no\end{align}
where $N$ is the random measure associated to $W_{S_t}$, i.e.,
\begin{align}
N(t,\Gamma):=\sum_{s\in(0,t]}1_{\Gamma}(W_{S_s}-W_{S_{s-}}),\ \ \Gamma\in\sB(\mR^d).\label{FF99}
\end{align}
\item[{\rm (ii)}]  For any $p>1$, there exists a constant $C>0$ such that for all $t>0, x\in \mR^d$ and $f\in\mathcal B_b(\mR^d)$,
\begin{align}
|\nabla P_t f(x)|\leq C(1+|x|^m)  (P_t|f|^p)^{1/p}(x)\left(\mE\left\{S_{t\land 1}^{-p/(2(p-1))}\vee 1\right\}\right)^{(p-1)/p}.\label{FG88}
\end{align}
Here $m$ is from (\ref{Ass1}).
Consequently, if $\liminf_{u\to\infty} B(u) u^{-\ff\aa 2}>0$  for some $\aa\in (0,2)$,
 then for any $p>1$, there exists a constant $C>0$ such that for all $t>0$, $x\in\mR^d$ and $f\in\mathcal B_b(\mR^d)$,
\begin{align}
|\nabla P_t f(x)|\leq \ff{C(1+|x|^m)}{(t\wedge 1)^{\frac{1}{\alpha}} }(P_t|f|^p)^{1/p}(x).\label{Ess}
\end{align}
\end{enumerate}
\et

\paragraph{Remark 1.1.} (1) We take the stopping time in \eqref{EE11} to make the stochastic integral $\int_0^{t\land \tau}\<\sigma^{-1}_s(X_{s-})\nabla_\v X_{s-},\dif W_{S_s}\>$ well-defined and square integrable. Since $S_{t\land\tau}\ge R$ for $\tau\le t$, we have $S^{-\ff 1 2}_{t\land\tau} \le S_t^{-\ff 1 2} +   R^{-\ff 1 2}$, so that the condition $\mE S_{t}^{-1/2}<\infty$ is equivalent to
$\mE S_{t\land\tau}^{-1/2}<\infty$ by the increasing property of $t\mapsto S_t$.

(2) According to known derivative estimates of the $\aa$-stable process, the gradient estimate (\ref{Ess})   is sharp in short-time for $S_t$ being the $\aa$-stable subordinator (i.e.  $B(r)=   r^{\ff\aa 2}$). Moreover, (\ref{Ess}) recovers (\ref{Z}) in the additive noise case by taking $m=0$. See e.g. \cite[Theorem 3.1(d)]{DS} for a condition on the symbol which is equivalent to
$\liminf_{u\to\infty} B(u) u^{-\ff\aa 2}>0$ used in Theorem 1.1(ii).

(3) When the noise is additive, the dimension-free Harnack inequality initiated in \cite{W97} was established in \cite{WW}
by using regularization approximations and coupling (see also \cite{WZ}). 
According to \cite[Propostion 1.3.1]{WB} (see also \cite{ATW09, GW}), this type Harnack inequality of $P_t$ follows from the entropy type gradient estimate
\beg{equation}\label{R*} |\nn P_t f|(x)\le \dd \big\{P_t(f\log f)-P_tf)\log P_t f\big\}(x) +C_t(x,\dd) P_t f(x),\ \ \dd>0, x\in \mR^d\end{equation} for positive $f\in \B_b(\R^d).$ However, this type gradient estimate can not be derived from the   formula \eqref{EE11}. On the other hand, in the multiplicative Gaussian noise case, the dimension-free Harnack inequality has been derived in \cite{W11b} by using coupling with an additional unbounded drift. It might be possible to modify this argument for the present setting.

(4) As a counterpart of the 
Bismut-Elworthy-Li 
type formula and dimension-free Harnack inequality, Driver type integration by parts formula and shift-Harnack inequality have been investigated in \cite{W14} for SDEs/SPDEs driven by Gaussian noise. We believe that under reasonable conditions the integration by parts formula can also be derived for the present setting by using finite-jump approximations introduced in this paper. To this end, for a fixed vector $v$ one needs to construct $h$ such that, instead of \eqref{ET2} below,  the Malliavin derivative $D_h X_t= v$. However, the construction of $h$ would be  more complicated than that for \eqref{ET2}, as already explained in \cite{W14} for the diffusion case. The study of shift-Harnack inequality in the multiplicative noise case would be  more difficult since the required backward coupling in \cite[Theorem 2.1]{W14} is not yet available even in the multiplicative Gaussian noise case. We note that in the additive noise case the shift-Harnack inequality has been derived in \cite{DS}. 
Moreover, some extensions for irregular drifts and other type multiplicative noises can be found in \cite{CSZ} and \cite{WXZ}.

\

The remainder of the paper is organized as follows. In Section 2 we briefly recall the integration by parts formula in the Malliavin calculus and present some lemmas on finite-jump approximations. These are then used in Section 3 to establish a derivative formula for $P_t^\ell$. Finally, in Section 4 we present explicit gradient estimates of $P_t^\ell$ and a complete proof of the above main result.

\section{Preliminaries}

Let $(\mW,\mH,\mu^\mW)$ be the classical Wiener space, i.e., $\mW$ is
the space of all continuous functions $\omega:\mR_+\to\mR^d$ with $\omega_0=0$;
$\mH$ is the Cameron-Martin space consisting of all absolutely continuous functions $h\in\mW$ with
$$
\|h\|^2_\mH:=\int^\infty_0 |\dot h_s|^2\dif s<+\infty,
$$
where $\dot h_s$ denotes the generalized derivative of $h_s$ with respect to $s$;
and $\mu^\mW$ is the Wiener measure so that the coordinate process
$$
W_t(\omega):=\omega_t
$$
is a $d$-dimensional standard Brownian motion.

Below we recall some basic notions about the Malliavin calculus (cf. \cite{Ma, Nu}). Let $\mU$ be a real separable Hilbert space.
Let $\sC(\mU)$ be the class of all $\mU$-valued smooth cylindrical functionals on $\Omega$ with the form:
$$
F=\sum_{i=1}^mf_i(W(h_1),\cdots,W(h_n)) u_i,
$$
where $m,n\ge 1, f_i\in C^\infty_0(\mR^n)$, $u_i\in\mU$, $h_1,\cdots, h_n\in\mH$ and
$$
W(h):=\int^\infty_0 \dot h_s\dif W_s.
$$
The Malliavin derivative of $F$ is defined by
$$
DF:=\sum_{i=1}^m\sum_{j=1}^n(\p_j f_i)(W(h_1),\cdots,W(h_n))u_i\otimes h_j\in\mU\otimes\mH.
$$
It is well known that the operator $(D,\sC(\mU))$ is closable from $L^2(\mW;\mU)$ to $L^2(\mW;\mU\otimes\mH)$ (cf. \cite[p.26, Proposition 1.2.1]{Nu}).
The closure is denoted by $(D,\sD(D))$. The dual operator $\delta$ of $D$ (also called divergence operator) is defined by
$$
\mE \<DF,U\>_\mH=\mE (F\delta( U)),\ \ F\in\sD(D),\ \ U\in \sD(\delta).
$$
Notice that the following divergence formula holds: for any $F\in\sD(D)$ and $U\in\sD(\dd)$,
\begin{align}
\delta(F U)=F\delta U-\<DF,U\>_\mH.\label{For}
\end{align}
Below we also write
$$
D_U F:=\<DF, U\>_\mH
$$
and
$$
\sF^\mW_t:=\sigma\{W_s: s\leq t\}.
$$
If $U_t(\omega)=\int^t_0u_s(\omega)\dif s$, where $u_s$ is an $\sF^\mW_t$-adapted
process on $\mR^d$ with $\mE\int^\infty_0|u_s|^2\dif s<\infty$, then $u\in\sD(\delta)$ and (cf. \cite{Nu})
$$
\delta (U)=\int^\infty_0\<u_s,\dif W_s\>.
$$

Next, for  $\beta: [0,\infty)\to [0,\infty)$ being an  absolutely continuous increasing function with $\beta_0=0$ and locally bounded derivative $\dot\bb_t$, set
\begin{align}
W^\beta_t:=\int^t_0\dot\beta_s\dif W_s,\ \ t\ge 0.\label{WB}
\end{align}
It is easy to see that $t\mapsto W^\beta_t$ is a process with independent increments and
\begin{align}
\lambda^\beta_t:=\mE |W^\beta_t|^2=\int^t_0|\dot\beta_s|^2\dif s.\label{LB}\end{align}
Let $\mS$ be the space of all purely jump  c\`adl\`ag and increasing functions, i.e. the path space of $S_t$. For any $\ell\in\mS$ and $s>0$, we shall denote
$$
\Delta\ell_s:=\ell_s-{\ell_{s-}}.
$$

We recall the following Burkholder's inequality (cf. \cite[Theorem 2.3]{Zh1}).
\bl\label{Le0}
Assume that $\xi_t$ is an $\sF^\mW_{\ell_t}$-adapted c\`adl\`ag $\mR^d$-valued process
and satisfies that for some $p>0$,
\begin{align}
\mE\left(\int^T_0|\xi_{s-}|^2\dif \lambda^\beta_{\ell_s}\right)^{\frac{p}{2}}<+\infty,\ \ \forall T\geq 0.\label{ET8}
\end{align}
Then there exists a constant $C_p>0$ such that for all $T\geq 0$,
\begin{align}
\mE\left(\sup_{t\in[0,T]}\left|\int^t_0\xi_{s-}\dif W^\beta_{\ell_s}\right|^p\right)\leq
C_p\mE\left(\int^T_0|\xi_{s-}|^2\dif \lambda^\beta_{\ell_s}\right)^{\frac{p}{2}}.\label{EU2}
\end{align}
\el

Finally, for fixed $\ell\in\mS$ and $\vv>0$,  define
\begin{align}
\ell^\eps_t:=\sum_{s \leq t} \Delta \ell_s 1_{\{\Delta \ell_s \geq \eps\}}, \ \ \ t\geq 0.\label{Ap}
\end{align}

We present below a key lemma about the approximation $\ell^\eps_t$ for later use.
\bl\label{Le1}
Let $\xi^\eps_t$ be an $\sF^\mW_{\ell^\eps_t}$-adapted c\`adl\`ag $\mR^d$-valued process such that
\begin{align}
\sup_{t\in[0,T]}\sup_{\eps\in(0,1)}\mE|\xi_t^\eps|^2<\infty,\ T>0, \label{TR1}
\end{align}
and for any $t>0$ and $\eps\in(0,1)$,
\begin{align}
\lim_{s\uparrow t}\mE|\xi_s^\eps-\xi_{t-}^\eps|^2=0. \label{TR11}
\end{align}
Then for any $T>0$, we have
$$
\lim_{\eps\downarrow 0}\sup_{t\in[0,T)}\mE\left|\int^t_0\xi^\eps_{s-}\dif W^\beta_{\ell^\eps_s}-\int^t_0\xi^\eps_{s-}\dif W^\beta_{\ell_s}\right|^2=0,
$$
where $W^\beta_t$ is defined by (\ref{WB}).
\el
\begin{proof}
Without loss of generality, we assume $T=1$. For $t\in[0,1)$,
set
$$
t_n:=[nt]/n,\ \ t_n^+:=([nt]+1)/n,
$$
where $[a]$ denotes the integer part of a real number $a$.
Fix $\eps\in(0,1)$. Since $\sF^\mW_{\ell^\eps_t}\subset\sF^\mW_{\ell_t}$, by (\ref{EU2}), (\ref{TR1}) and (\ref{TR11}), we have
\begin{align}
&\mE\left|\int^t_0\xi^\eps_{s-}\dif W^\beta_{\ell^\eps_s}-\int^t_0\xi^\eps_{s-}\dif W^\beta_{\ell_s}\right|^2
=\lim_{n\to\infty}\mE\left|\int^{t}_0\xi^\eps_{s_n-}\dif W^\beta_{\ell^\eps_s}
-\int^{t}_0\xi^\eps_{s_n-}\dif W^\beta_{\ell_s}\right|^2\no\\
&\qquad=\lim_{n\to\infty}\mE\left|\sum_{k=0}^{[nt]}\xi^\eps_{k/n}\left(W^\beta_{\ell^\eps_{(k+1)/n}}
-W^\beta_{\ell^\eps_{k/n}}-W^\beta_{\ell_{(k+1)/n}}+W^\beta_{\ell_{k/n}}\right)\right|^2\no\\
&\qquad=\lim_{n\to\infty}\sum_{k=0}^{[nt]}\mE|\xi^\eps_{k/n}|^2\mE\left(W^\beta_{\ell^\eps_{(k+1)/n}}
-W^\beta_{\ell^\eps_{k/n}}-W^\beta_{\ell_{(k+1)/n}}+W^\beta_{\ell_{k/n}}\right)^2\no\\
&\qquad\leq C\lim_{n\to\infty}\sum_{k=0}^{n-1}\mE\left(W^\beta_{\ell^\eps_{(k+1)/n}}
-W^\beta_{\ell^\eps_{k/n}}-W^\beta_{\ell_{(k+1)/n}}+W^\beta_{\ell_{k/n}}\right)^2,\label{EY1}
\end{align}
where we have used the independence of $\xi^\eps_{k/n}$ and $W^\beta_{\ell^\eps_{(k+1)/n}}
-W^\beta_{\ell^\eps_{k/n}}-W^\beta_{\ell_{(k+1)/n}}+W^\beta_{\ell_{k/n}}$.

Notice that if $\ell_{k/n}\geq\ell^\eps_{(k+1)/n}$, then
\begin{align*}
&\mE\left(W^\beta_{\ell^\eps_{(k+1)/n}}-W^\beta_{\ell^\eps_{k/n}}-W^\beta_{\ell_{(k+1)/n}}+W^\beta_{\ell_{k/n}}\right)^2
=\mE\left(\int^{\ell^\eps_{(k+1)/n}}_{\ell^\eps_{k/n}}\dot\beta_s\dif W_s-\int^{\ell_{(k+1)/n}}_{\ell_{k/n}}\dot\beta_s\dif W_s\right)^2\\
&\quad=\int^{\ell^\eps_{(k+1)/n}}_{\ell^\eps_{k/n}}|\dot\beta_s|^2\dif s+\int^{\ell_{(k+1)/n}}_{\ell_{k/n}}|\dot\beta_s|^2\dif s
=\lambda^\beta_{\ell^\eps_{(k+1)/n}}-\lambda^\beta_{\ell^\eps_{k/n}}+\lambda^\beta_{\ell_{(k+1)/n}}-\lambda^\beta_{\ell_{k/n}},
\end{align*}
and if $\ell_{k/n}<\ell^\eps_{(k+1)/n}$, then
\begin{align*}
&\mE\left(W^\beta_{\ell^\eps_{(k+1)/n}}-W^\beta_{\ell^\eps_{k/n}}-W^\beta_{\ell_{(k+1)/n}}+W^\beta_{\ell_{k/n}}\right)^2
=\mE\left(\int^{\ell_{k/n}}_{\ell^\eps_{k/n}}\dot\beta_s\dif W_s-\int^{\ell_{(k+1)/n}}_{\ell^\eps_{(k+1)/n}}\dot\beta_s\dif W_s\right)^2\\
&\quad=\int^{\ell_{k/n}}_{\ell^\eps_{k/n}}|\dot\beta_s|^2\dif s+\int^{\ell_{(k+1)/n}}_{\ell^\eps_{(k+1)/n}}|\dot\beta_s|^2\dif s
=\lambda^\beta_{\ell_{k/n}}-\lambda^\beta_{\ell^\eps_{k/n}}+\lambda^\beta_{\ell_{(k+1)/n}}-\lambda^\beta_{\ell^\eps_{(k+1)/n}}.
\end{align*}
Hence,
\begin{align}
\sum_{k=0}^{n-1}\mE\left(W^\beta_{\ell^\eps_{(k+1)/n}}
-W^\beta_{\ell^\eps_{k/n}}-W^\beta_{\ell_{(k+1)/n}}+W^\beta_{\ell_{k/n}}\right)^2=I_1(n,\eps)+I_2(n,\eps),\label{EY2}
\end{align}
where
\begin{align*}
I_1(n,\eps)&:=\sum_{k=0}^{n-1}\left(\lambda^\beta_{\ell^\eps_{(k+1)/n}}-\lambda^\beta_{\ell^\eps_{k/n}}+\lambda^\beta_{\ell_{(k+1)/n}}-\lambda^\beta_{\ell_{k/n}}\right)1_{\ell_{k/n}\geq\ell^\eps_{(k+1)/n}}\\
&=\int^1_0\left(\frac{\lambda^\beta_{\ell^\eps_{s^+_n}}-\lambda^\beta_{\ell^\eps_{s_n}}}{\ell_{s^+_n}-\ell_{s_n}}
+\frac{\lambda^\beta_{\ell_{s^+_n}}-\lambda^\beta_{\ell_{s_n}}}{\ell_{s^+_n}-\ell_{s_n}}\right)1_{\ell_{s_n}\geq\ell^\eps_{s^+_n}}\dif\ell_s,\\
I_2(n,\eps)&:=\sum_{k=0}^{n-1}\left(\lambda^\beta_{\ell_{k/n}}-\lambda^\beta_{\ell^\eps_{k/n}}+\lambda^\beta_{\ell_{(k+1)/n}}-\lambda^\beta_{\ell^\eps_{(k+1)/n}}\right)
1_{\{\ell_{k/n}<\ell^\eps_{(k+1)/n}\}}\\
&=\int^1_0\left(\frac{\lambda^\beta_{\ell_{s_n}}-\lambda^\beta_{\ell^\eps_{s_n}}}{\ell_{s^+_n}-\ell_{s_n}}
+\frac{\lambda^\beta_{\ell_{s^+_n}}-\lambda^\beta_{\ell^\eps_{s^+_n}}}{\ell_{s^+_n}-\ell_{s_n}}\right)1_{\{\ell_{s_n}<\ell^\eps_{s^+_n}\}}\dif\ell_s.
\end{align*}
Noticing that as $n\to\infty$,
$$
\ell_{s_n}\uparrow\ell_{s-},\ \ \ell_{s^+_n}\downarrow\ell_s,\ \ \ell^\eps_{s_n}\uparrow\ell^\eps_{s-},\ \ \ell^\eps_{s^+_n}\downarrow\ell^\eps_s,
$$
and as $\eps\downarrow 0$,
$$
\ell^\eps_s\uparrow\ell_s,\ \ \ell_{s-}^\vv\uparrow \ell_{s-},
$$
by the dominated convergence theorem, we have
\begin{align}
\lim_{\eps\downarrow 0}\varlimsup_{n\to\infty}I_1(n,\eps)&=\lim_{\eps\downarrow 0}\int^1_0\left(\frac{\lambda^\beta_{\ell^\eps_{s}}-\lambda^\beta_{\ell^\eps_{s-}}}{\ell_{s}-\ell_{s-}}
+\frac{\lambda^\beta_{\ell_{s}}-\lambda^\beta_{\ell_{s-}}}{\ell_{s}-\ell_{s-}}\right)1_{\{\ell_{s-}\geq\ell^\eps_{s}\}}\dif\ell_s\no\\
&=2\int^1_0\left(\frac{\lambda^\beta_{\ell_{s}}-\lambda^\beta_{\ell_{s-}}}{\ell_{s}-\ell_{s-}}\right)1_{\{\ell_{s-}\geq\ell_{s}\}}\dif\ell_s=0,\label{EY3}
\end{align}
and
\begin{align}
\lim_{\eps\downarrow 0}\varlimsup_{n\to\infty}I_2(n,\eps)&\leq\lim_{\eps\downarrow 0}\int^1_0\left(\frac{\lambda^\beta_{\ell_{s-}}-\lambda^\beta_{\ell^\eps_{s-}}}{\ell_{s}-\ell_{s-}}
+\frac{\lambda^\beta_{\ell_{s}}-\lambda^\beta_{\ell^\eps_{s}}}{\ell_{s}-\ell_{s-}}\right)1_{\{\ell_{s-}\leq \ell^\vv_s\}} \dif\ell_s =0.\label{EY4}
\end{align}
Combining (\ref{EY1})-(\ref{EY4}), we obtain the desired limit.
\end{proof}

\section{Derivative formula for $P_t^\ell$}

In this section, we fix an $\ell\in\mS$ and consider the equation
\beq\label{3.0}
X^\ell_t(x)=x+\int^t_0b_s(X^\ell_s(x))\dif s+\int^t_0\sigma_s(X^\ell_{s-}(x))\dif W_{\ell_s}
\end{equation}  for $b$ and $\si$ satisfying {\bf (A)}.
Under (\ref{Ass0}), it is well-known that $\{X^\ell_t(\cdot),t\geq 0\}$ forms a $C^1$-stochastic flow (cf. \cite[p.305, Theorem 39]{Pr}).
Let $\nabla X^\ell_t:=(\p_j (X^\ell_t)_i(x))_{ij}$ be the derivative matrix. Then
$$\nn_\v X^\ell_t:= (\nn_\v (X^\ell_t)_i)_{1\leq i\leq d}:= (\nn X^\ell_t)\v,\ \ \v\in\mR^d,$$
and
\beg{equation}\label{TR2}
\nabla_\v X^\ell_t=\v+\int^t_0\nabla_{\nabla_\v X^\ell_s}  b_s(X^\ell_s)\dif s+\int^t_0\nabla_{\nabla_\v X^\ell_{s-}} \sigma_s(X^\ell_{s-})\dif W_{\ell_s}.
\end{equation}

We first prepare the following lemma.
\bl\label{Le2} Assume  {\bf (A)}. Let $\aa_t:= \sup_{s\in [0, t]} |b_s(0)|$ and $\gg_t:=  \sup_{s\in [0, t]} \|\si_s(0)\|$. Then,
for any $p\geq 2$ there exists a constant $c(p)>0$ such that for all $t\geq 0$ and $x\in\mR^d$,
\begin{align}
\mE\left(\sup_{s\in [0,t]} |X_s^\ell(x)|^p\right) \leq c(p) \Big(|x|+\aa_t t+ \gg_t \sqrt{\ell_t}\Big)^p \exp\left\{{c(p)(\|\nn b\|_t^pt^p+\|\nn \si\|_t^p\ell_t^{\ff p2})}\right\},\label{L1.4.1}
\end{align}
and
\begin{align}
\mE\left(\sup_{s\in [0,t]} \|\nn X_s^\ell(x)\|^p\right) \leq c(p)   \exp\left\{c(p)(\|\nn b\|_t^pt^p+\|\nn \si\|_t^p\ell_t^{\ff p 2})\right\}.\label{L1.4.2}\end{align}
Moreover, letting $X_t^{\ell^\vv}$ solving \eqref{3.0} for $\ell^\vv$ in place of $\ell$,  for any $T>0$ and $x\in\mR^d$ we have
\begin{align}
\lim_{\eps\downarrow 0}\sup_{t\in[0,T]}\mE|X^{\ell^\eps}_t(x)-X^{\ell}_t(x)|^2=0,\ \
\lim_{\eps\downarrow 0}\sup_{t\in[0,T]}\mE|\nabla X^{\ell^\eps}_t(x)-\nabla X^{\ell}_t(x)|^2=0.\label{L1.4.3}
\end{align}
\el
\begin{proof} (1) We simply use $X_s$ to denote $X_s^\ell(x).$
 Since
$$|X_s|^p\leq 3^{p-1}|x| + 3^{p-1} \bigg|\int_0^t b_s(X_s) \dif s\bigg|^p + 3^{p-1} \bigg|\int_0^t \si_s(X_{s-}) \dif W_{\ell_s}\bigg|^p$$ and
$$|b_s(x)|\leq \|\nn b\|_s |x|+\aa_s,\ \ \|\si_s(x)\|\leq \|\nn \si\|_s |x|+\gg_s, $$
by Lemma \ref{Le0} and H\"older's inequality, there exists a constant $c(p)>0$ such that
\begin{align*} \mE \left(\sup_{s\in [0,t]} |X_s|^p\right) \leq & c(p)\Big(|x|+ t\aa_t +\gg_t\sqrt{\ell_t}\Big)^p + c(p)\|\nn b\|_t^pt^{p-1}\int_0^t \mE\left(\sup_{r\in [0,s]} |X_r|^p\right) \dif s\\
&+ c(p)\|\nabla\si\|_t^p\ell_t^{\ff p 2-1}\int_0^t \mE\left(\sup_{r\in [0,s]} |X_r|^p\right) \dif \ell_s.\end{align*}
By Gronwall's lemma, this implies (\ref{L1.4.1}).

Similarly,
noticing that \eqref{TR2} implies
$$
\|\nabla X_t\|^p\leq 3^{p-1} +3^{p-1} \|\nn b\|_t^p \bigg|\int^t_0\|\nabla X_s\|\dif s\bigg|^p
+3^{p-1} \bigg\|\int^t_0\nabla\sigma_s(X^\ell_{s-})\nabla X^\ell_{s-}\dif W_{\ell_s}\bigg\|^p,
$$
by Lemma \ref{Le0} again and H\"older's inequality, there exists a constant $c(p)>0$ such that
\begin{align*}
\mE\left(\sup_{s\in [0,t]}\|\nabla X_t\|^p\right)&\leq  c(p)+ c(p)\|\nn b\|_t^pt^{p-1}\int_0^t \mE\left(\sup_{r\in [0,s]} \|\nn X_r\|^p\right) \dif s\\
&+ c(p)\|\nn \si\|_t^p\ell_t^{\ff p 2-1}\int_0^t \mE\left(\sup_{r\in [0,s]} \|\nn X_r\|^p\right) \dif \ell_s,
\end{align*}
which implies (\ref{L1.4.2}) by Gronwall's inequality again.
\\
\\
(2) Let  $Y^\eps_t:=X^{\ell^\eps}_t-X^{\ell}_t.$ Then for all $t\leq T$,
$$
\mE|Y^\eps_t|^2\leq 3 \|\nabla b\|^2_t T\int^t_0\mE|Y^\eps_s|^2\dif s
+3\|\nabla \sigma\|^2_t\int^t_0\mE|Y^\eps_s|^2\dif \ell_s+3R^\eps_t,
$$
where
$$
R^\eps_t:=\mE\left|\int^t_0\sigma_s(X^\eps_{s-})\dif W_{\ell^\eps_s}-\int^t_0\sigma_s(X^\eps_{s-})\dif W_{\ell_s}\right|^2.
$$
By Gronwall's lemma and Lemma \ref{Le1}, we obtain
$$
\lim_{\eps\downarrow 0}\sup_{t\in[0,T]}\mE|Y^\eps_t|^2\leq \Big(3\,\e^{3\|\nn b\|^2_T T^2+3\|\nn \sigma\|^2_T\ell_T}\Big)\,\lim_{\eps\downarrow 0}\sup_{t\in[0,T]}R^\eps_t=0.
$$
Then the first part in \eqref{L1.4.3} holds. The proof for the second part is completely similar.
\end{proof}

We now prove the following derivative formula:

\bt\label{Main} Assume {\bf (A)}.  Let $\beta$ be an increasing $C^1$-function with $\beta_0=0$.
For any $\mathrm{v} \in \mR^d$ and $f\in \mathcal B_b(\mR^d)$, we have for all $t>0$ with $\beta_{\ell_t}>0$,
\begin{align}
\nabla_{\mathrm{v}} P_t^\ell f=\frac{1}{\beta_{\ell_t}}\mE\Bigg[& f(X^{\ell}_t)\Bigg(\int^t_0
\<\sigma^{-1}_s(X^\ell_{s-})\nabla_\v X^\ell_{s-},\dif W^\beta_{\ell_s}\>\no\\
&-\int^t_0\mathrm{Tr}\big(\sigma^{-1}_s\nabla_{\nabla_\v X^\ell_{s-}}\sigma_s\big) (X^\ell_{s-}) \dif\beta_{\ell_s}\label{EE1}\\
& +\sum_{s\in[0, t]}\<\sigma^{-1}_s(X^\ell_{s-})\nabla_{\nabla_\v X^\ell_{s-}}\sigma_s(X^\ell_{s-})\Delta W^\beta_{\ell_{s}},\Delta W_{\ell_{s}}\>\Bigg)\Bigg],\no
\end{align}
where $W^\beta$ is defined by (\ref{WB}), and $\Delta W^\beta_{\ell_{s}}:=W^\beta_{\ell_{s}}-W^\beta_{\ell_{s-}}$.
\et
\begin{proof}  (i)
We first show that
it suffices to prove \eqref{EE1} for  $f\in C_0^\infty(\R^d).$  By {\bf (A)} and Lemma \ref{Le2},
\beg{equation*} \beg{split} \eta:= &\int^t_0
\<\sigma^{-1}_s(X^\ell_{s-})\nabla_\v X^\ell_{s-},\dif W^\beta_{\ell_s}\>
 -\int^t_0\mathrm{Tr}\big(\sigma^{-1}_s\nabla_{\nabla_\v X^\ell_{s-}}\sigma_s\big) (X^\ell_{s-}) \dif\beta_{\ell_s}\\
& +\sum_{s\in[0, t]}\<\sigma^{-1}_s(X^\ell_{s-})\nabla_{\nabla_\v X^\ell_{s-}}\sigma_s(X^\ell_{s-})\Delta W^\beta_{\ell_{s}},\Delta W_{\ell_{s}}\>\end{split}\end{equation*}
has finite moment. If \eqref{EE1} holds for $f\in C_0^\infty(\R^d)$ then
 \beq\label{EE1'} \nn_\v P_t f(x)= \ff 1 {\bb_{\ell_t}}\E\big\{f(X_t^\ell(x)\eta\big\},\ \ f\in C_0^\infty(\R^d).\end{equation}
 Since $\bb_{\ell_t}>0$ implies $\ell_t>0$, there exists a jump point $t_0\in (0,t].$ Obviously,
$$X_{t_0}^\ell(x)= X_{t_0-}^\ell(x) +\si_{t_0}(X_{t_0-}^\ell(x))\DD W_{\ell_{t_0}}.$$
Since $\DD W_{\ell_{t_0}}$ is a non-trivial normal random variable independent of $X_{t_0-}^\ell(x)$ and since $\si_{t_0}$ is invertible,
$X_{t_0}^\ell(x)$ has a density
$$p_0(x,y)= \E\Bigg\{\ff {\exp\big[-\ff 1 {2\DD\ell_{t_0}}\big|\si_{t_0}^{-1}(X_{t_0-}^\ell(x))(y-X_{t_0-}(x))\big|^2\big]}
{(2\pi\DD\ell_{t_0})^{d/2} \ss{\text{det}\big((\si_{t_0}\si_{t_0}^*)(X_{t_0-}^\ell(x))}\big)} \Bigg\}.$$
From {\bf (A)} and Lemma \ref{Le2}
we see that $x\mapsto p_0(x,y)$ is differentiable and
\beq\label{**} x\mapsto \int_{\R^d} |\nn p(\cdot,y)|(x)\d y\ \text{is\ locally\ bounded}.\end{equation}  Letting $X_{t_0,s}^\ell(z)$ solve the equation
$$ X_{t_0,s}^\ell(z)=z +\int_{t_0}^s b_r(X_{t_0,r}^\ell (z))\d r +\int_{t_0}^s \si_s(X_{t_0,r-}^\ell (z)) \d W_{\ell_r},\ \ s\in [t_0,t],$$ and denoting
$P^\ell_{t_0,t}f(z)= \E f(X_{t_0,t}^\ell(z))$,
by the Markov property we have
$$
P_t^\ell f(x)= P^\ell _{t_0}(P^\ell_{t_0,t}f)(x)=\int_{\R^d}(P^\ell_{t_0,t}f)(y)p_0(x,y)\d y.
$$
Combining this with  \eqref{**}, we may apply  the dominated convergence theorem to derive
\beq\label{0*}
\nn_\v P_t^\ell f(x)= \int_{\R^d} \{\nn_\v p_0(\cdot,y)(x) \} P^\ell_{t_0,t} f(y)\d y,\ \ f\in \B_b(\R^d).
\end{equation}
Moreover, \eqref{**} implies that
$$
\mu(A):= \int_{\R^d} \{\nn_\v p(\cdot,y)(x)\}P_{t_0,t}1_A(y) \d y- \ff 1 {\bb_{\ell_t}}\E\big\{1_A(X_t^\ell(x))\eta\big\}
$$
gives a finite signed measure on $\R^d$.  Since \eqref{EE1'} and  \eqref{0*} imply $\int_{\R^d} f\d\mu=0$ for all $f\in C_0^\infty(\R^d)$, we have $\mu=0$ so that
$$0=\int_{\R^d} f\d\mu = \nn_\v P_t^\ell f(x)-\ff 1 {\bb_{\ell_t}}\E\big\{f(X_t^\ell(x))\eta\big\},\  \ f\in \B_b(\R^d).$$
Thus, \eqref{EE1} holds for  $f\in \B_b(\R^d).$

(ii) Let us  establish formula (\ref{EE1}) for $\ell\in\mS$ with {\it finite many jumps} on any finite time interval.
Let $\cJ_t$ solve the following matrix-valued ODE
$$
\cJ_t=I+\int^t_0\nabla b_s(X_s)\cJ_s\dif s.
$$
By (\ref{TR2}) and the variation of constant formula, it is easy to see that
\begin{align}
\nabla_\v X_t=\cJ_t\v+\cJ_t\int^t_0\cJ^{-1}_s\nabla_{\nabla_\v X_{s-}}\sigma_s(X_{s-})\dif W_{\ell_s},\label{ET1}
\end{align}
where $\cJ^{-1}_s$ denotes the inverse matrix of $\cJ_s$.
Let $0=t_0<t_1<t_2<\cdots<t_n<\cdots$ be the jump times of $\ell$. Fix $\v\in\mR^d$ and set $h_0=0$. Define
$h\in\mH$ recursively as follows: if $s\in (\ell_{t_{k-1}},\ell_{t_k}]$ for some $k\geq 1$, then
\begin{align}
h_s:=h_{\ell_{t_{k-1}}}+(\beta_{s}-\beta_{\ell_{t_{k-1}}})\sigma^{-1}_{t_k}(X_{t_k-})\(\nabla_\v X_{t_k-}+\nabla_{\nabla_\v X_{t_k-}}\sigma_{t_k}(X_{t_k-})\Delta W_{\ell_{t_k}}\).\label{ET4}
\end{align}
For this $h$, recalling that $D_hX_t:=\<DX_t,h\>_\mH$, it is standard to prove that (cf. \cite{Nu})
$$
D_hX_t=\int^t_0\nabla_{D_h X_s} b_s(X_s)\dif s+\int^t_0\nabla_{D_h X_{s-}}\sigma_s(X_{s-})\dif W_{\ell_s}  +\int^t_0\sigma_s(X_{s-})\dif h_{\ell_s}.
$$
As above, one can write
\begin{align}
D_hX_t&=\cJ_t\int^t_0\cJ^{-1}_s\nabla_{D_h X_{s-}}\sigma_s(X_{s-})\dif W_{\ell_s}  +\cJ_t\int^t_0\cJ^{-1}_s\sigma_s(X_{s-})\dif h_{\ell_s}\no\\
&=\cJ_t\sum_{k=1}^m\cJ^{-1}_{t_k}\left(\nabla_{D_h X_{t_k-}}\sigma_{t_k}(X_{t_k-})\Delta W_{\ell_{t_k}}
+\sigma_{t_k}(X_{t_k-})\Delta h_{\ell_{t_k}}\right),\label{ET3}
\end{align}
where $t\in[t_m,t_{m+1})$, $\Delta W_{\ell_{t_k}}:=W_{\ell_{t_k}}-W_{\ell_{t_k-}}$ and $\Delta h_{\ell_{t_k}}:=h_{\ell_{t_k}}-h_{\ell_{t_k-}}$.
\\
\\
Now, we use the induction method to prove
\begin{align}
D_hX_t=\beta_{\ell_t}\nabla_\v X_t,\ \ t\geq 0.\label{ET2}
\end{align}
First of all, by definition we have
$$
D_hX_t=0,\ \ \beta_{\ell_t}=0,\ \ \forall t\in[0,t_1),
$$
and so (\ref{ET2}) holds for all $t\in[0,t_1)$.
Supposing that (\ref{ET2}) holds for all $t\in[0,t_m)$ for some $m\ge 1$, we aim to prove it for $t\in [t_m, t_{m+1}).$  Observing that
\begin{align}
\ell_{t_k-}=\ell_{t_{k-1}},\label{TR3}
\end{align}
by (\ref{ET4}) we have
$$
\Delta h_{\ell_{t_k}}=\Delta\beta_{\ell_{t_k}}\sigma^{-1}_{t_k}(X_{t_k-})\(\nabla_\v X_{t_k-}+\nabla_{\nabla_\v X_{t_k-}}\sigma_{t_k}(X_{t_k-})\Delta W_{\ell_{t_k}}\).
$$
Substituting this and induction hypothesis (\ref{ET2}) into  (\ref{ET3}), we obtain that for $t\in[t_m,t_{m+1})$,
\begin{equation}\label{GE}\begin{split}
D_hX_t&=\cJ_t\sum_{k=1}^m\cJ^{-1}_{t_k}\bigg\{\nabla_{\beta_{\ell_{t_k-}}\nabla_\v X_{t_k-}}\sigma_{t_k}(X_{t_k-})\Delta W_{\ell_{t_k}}\\
&\qquad +\Delta\beta_{\ell_{t_k}}\(\nabla_\v X_{t_k-}+\nabla_{\nabla_\v X_{t_k-}}\sigma_{t_k}(X_{t_k-})\Delta W_{\ell_{t_k}}\)\bigg\} \\
&=\cJ_t\sum_{k=1}^m\cJ^{-1}_{t_k}\bigg\{\beta_{\ell_{t_k}}\nabla_{\nabla_\v X_{t_k-}}\sigma_{t_k}(X_{t_k-})\Delta W_{\ell_{t_k}}
+\Delta\beta_{\ell_{t_k}}\nabla_\v X_{t_k-}\bigg\}  \\
&=\cJ_t\bigg(\int^t_0\cJ^{-1}_s\beta_{\ell_{s}}\nabla_{\nabla_\v X_{s-}}\sigma_s(X_{s-})\dif W_{\ell_{s}}
+\int^t_0\cJ^{-1}_s\nabla_\v X_{s-}\dif\beta_{\ell_s}\bigg),
\end{split}\end{equation} where in the second step we have used $ \Delta\beta_{\ell_{t_k}}=\beta_{\ell_{t_k}}-\beta_{\ell_{t_k-}}.$
On the other hand, since $s\mapsto\cJ_s$ is continuous, by (\ref{ET1}) we have
$$
\cJ^{-1}_s\nabla_\v X_{s-}=\v+\int_{[0,s)}\cJ^{-1}_r\nabla_{\nabla_\v X_{r-}}\sigma_r(X_{r-})\dif W_{\ell_r},
$$
and by Fubini's theorem,
\begin{align*}
\int^t_0\cJ^{-1}_s\nabla_\v X_{s-}\dif\beta_{\ell_s}
&=\beta_{\ell_t}\v+\int^t_0\!\!\!\int_{[0,s)}\cJ^{-1}_r\nabla_{\nabla_\v X_{r-}}\sigma_r(X_{r-})\dif W_{\ell_r}\dif\beta_{\ell_s}\\
&=\beta_{\ell_t}\v+\int^t_0(\beta_{\ell_t}-\beta_{\ell_r})\cJ^{-1}_r\nabla_{\nabla_\v X_{r-}}\sigma_r(X_{r-})\dif W_{\ell_r}.
\end{align*}
Substituting this into (\ref{GE}) and applying (\ref{ET1}), we arrive at
$$
D_hX_t=\cJ_t\left(\beta_{\ell_t}\v+\beta_{\ell_t}\int^t_0\cJ^{-1}_r\nabla_{\nabla_\v X_{r-}}\sigma_r(X_{r-})\dif W_{\ell_r}\right) =\beta_{\ell_t}\nabla_\v X_t.
$$
Thus, (\ref{ET2}) is proven.
\\
\\
Finally, assume that
$
\beta_{\ell_t}>0.
$
By (\ref{ET2}) and the integration by parts formula in the Malliavin calculus, we have
\begin{align}
\nabla_\v\mE f(X_t)&=\mE\(\nabla f(X_t)\nabla_\v X_t\)=\frac{1}{\beta_{\ell_t}}\mE\(\nabla f(X_t)D_h X_t\)\no\\
&=\frac{1}{\beta_{\ell_t}}\mE\(D_h(f(X_t))\)=\frac{1}{\beta_{\ell_t}}\mE\(f(X_t)\delta(h)\),\label{F0}
\end{align}
provided $h\in\sD(\dd)$.
Notice that $h$ defined in (\ref{ET4}) is {\it non-adapted} because of the term $\Delta W_{\ell_{t_k}}$. Neverethless, $\delta(h)$ can be explicitly calculated as follows:
Define
$$
h_k(s)=(\beta_{(s\lor \ell_{t_{k-1}})\land \ell_{t_k}}-\beta_{\ell_{t_{k-1}}})\sigma^{-1}_{t_k}(X_{t_k-})\(\nabla_\v X_{t_k-}+\nabla_{\nabla_\v X_{t_k-}}\sigma_{t_k}(X_{t_k-})\Delta W_{\ell_{t_k}}\),\ \ s\geq 0.
$$
Noting that
$$h_k= h_k^{(0)} +\sum_{j=1}^d h_k^{(j)}(\DD  W_{\ell_{t_k}})_j,$$
where
\begin{align*} & h_k^{(0)}(s) := (\beta_{(s\lor \ell_{t_{k-1}})\land \ell_{t_k}}-\beta_{\ell_{t_{k-1}}})\sigma^{-1}_{t_k}(X_{t_k-})\nabla_\v X_{t_k-},\\
&h_k^{(j)}(s):= (\beta_{(s\lor \ell_{t_{k-1}})\land \ell_{t_k}}-\beta_{\ell_{t_{k-1}}})\Big(\sigma^{-1}_{t_k}(X_{t_k-})\nabla_{\nabla_\v X_{t_k-}}\sigma_{t_k}(X_{t_k-})\Big)_{\cdot j}\end{align*} are adapted, we see from (\ref{Ass1}) and (\ref{L1.4.1}), (\ref{L1.4.2}) that $h\in\sD(\dd)$; and by (\ref{For}),
\begin{align*}
\delta(h)&=\sum_{k: t_k\leq t}\delta(h_k)=\sum_{k: t_k\leq t}\bigg(\dd(h_k^{(0)})+\sum_{j=1}^d \big\{\dd(h_k^{(j)})(\DD W_{\ell_{t_k}})_j-D_{h_k^{(j)}} (\DD W_{\ell_{t_k}})_j\big\}\bigg)\\
&=\sum_{k: t_k\leq t}\(\<\sigma^{-1}_{t_k}(X_{t_k-})\nabla_\v X_{t_k-},\Delta W^\beta_{\ell_{t_k}}\>+
\<\sigma^{-1}_{t_k}(X_{t_k-})\nabla_{\nabla_\v X_{t_k-}}\sigma_{t_k}(X_{t_k-})\Delta W^\beta_{\ell_{t_k}},\Delta W_{\ell_{t_k}}\>\)\\
&\quad -\sum_{k: t_k\leq t}(\beta_{\ell_{t_k}}-\beta_{\ell_{t_{k-1}}})\mathrm{Tr}\big(\sigma^{-1}_{t_k} \nabla_{\nabla_\v X_{t_k-}}\sigma_{t_k} \big)(X_{t_k-}) \\
&=\int^t_0\<\sigma^{-1}_s(X_{s-})\nabla_\v X_{s-},\dif W^\beta_{\ell_s}\>
- \int^t_0\mathrm{Tr}\big(\sigma^{-1}_s \nabla_{\nabla_\v X_{s-}}\sigma_s\big)(X_{s-})\dif\beta_{\ell_s}\no\\
&\qquad+\sum_{s\in [0, t]}\<\sigma^{-1}_s(X_{s-})\nabla_{\nabla_\v X_{s-}}\sigma_s(X_{s-})\Delta W^\beta_{\ell_{s}},\Delta W_{\ell_{s}}\>.
\end{align*}
By (\ref{F0}), this implies  (\ref{EE1}).
\\
\\
(iii) For general $\ell\in\mS$. Let $\ell^\eps$ be defined by (\ref{Ap}). It is well known that $\ell^\eps$ has finite many jumps on any finite interval.
By formula (\ref{EE1}) for $\ell^\eps$, we have
\begin{align}
\nabla_{\mathrm{v}} \mE f(X^{\ell^\eps}_t)=I_1(\eps)+I_2(\eps)+I_3(\eps),
\label{F000}\end{align}
where
\begin{align*}
I_1(\eps)&:=\frac{1}{\beta_{\ell^\eps_t}}\mE\left[f(X^{\ell^\eps}_t)
\int_0^t\sigma^{-1}_{s}(X^{\ell^\eps}_{s-})\left(\nabla_{\mathrm{v}} X^{\ell^\eps}_{s-}\right)\dif W^\beta_{\ell^\eps_{s}}\right],\no\\
I_2(\eps)&:=\frac{1}{\beta_{\ell^\eps_t}}\mE\left[f(X^{\ell^\eps}_t)\int_0^t\mathrm{Tr}\big(\sigma^{-1}_s
\nabla_{\nabla_\v X^{\ell^\eps}_{s-}}\sigma_s\big)(X^{\ell^\eps}_{s-})\dif\beta_{\ell^\eps_s}\right],\no\\
I_3(\eps)&:=\frac{1}{\beta_{\ell^\eps_t}}\mE\left[f(X^{\ell^\eps}_t)
\sum_{s\leq t}\<\sigma^{-1}_s(X^{\ell^\eps}_{s-})\nabla_{\nabla_\v X^{\ell^\eps}_{s-}}\sigma_s(X^{\ell^\eps}_{s-})\Delta W^\beta_{\ell^\eps_{s}},\Delta W_{\ell^\eps_{s}}\>\right].
\end{align*}
Set
$$
\xi^\eps_s=\sigma^{-1}_{s}(X^{\ell^\eps}_{s-}) \nabla_{\mathrm{v}} X^{\ell^\eps}_{s-},\ \
 \xi_s=\sigma^{-1}_{s}(X^{\ell}_{s-} )\nabla_{\mathrm{v}} X^{\ell}_{s-}.$$
By (\ref{Ass0}), (\ref{Ass1}) and Lemma \ref{Le2}, it is easy to see that
$$
\lim_{\eps\downarrow 0}\mE|\xi^\eps_s-\xi_s|^2=0.
$$
Hence, by Lemmas \ref{Le1}, \ref{Le2} and the dominated convergence theorem, we have
\begin{equation*}\begin{split}
&\lim_{\eps\downarrow 0}\mE\left|\int^t_0\xi^\eps_s\dif W^\beta_{\ell^\eps_s}-\int^t_0\xi_s\dif W^\beta_{\ell_s}\right|^2
\leq 2\lim_{\eps\downarrow 0}\mE\left|\int^t_0(\xi^\eps_s-\xi_s)\dif W^\beta_{\ell_s}\right|^2\\
&=2\lim_{\eps\downarrow 0}\int^t_0\mE|\xi^\eps_s-\xi_s|^2\dif \lambda^\beta_{\ell_s}=0,
\end{split}\end{equation*}
and so,
\begin{align}
\lim_{\vv\downarrow 0}I_1(\eps)=\frac{1}{\beta_{\ell_t}}\mE\left[f(X^{\ell}_t)
\int_0^t\sigma^{-1}_{s}(X^{\ell}_{s-})\left(\nabla_{\mathrm{v}} X^{\ell}_{s-}\right)\dif W^\beta_{\ell_{s}}\right].
\label{F1}\end{align}
Similarly, we have
\begin{align}
\lim_{\vv\downarrow 0}I_2(\eps)= \frac{1}{\beta_{\ell_t}}\mE\left[f(X^{\ell}_t)
\int^t_0\mathrm{Tr}\big(\sigma^{-1}_s \nabla_{\nabla_\v X^\ell_{s-}}\sigma_s\big)(X^\ell_{s-})\dif\beta_{\ell_s}\right].
\label{F2}\end{align}
To treat $I_3(\eps)$, we set
\begin{align}
A^\eps_s:= \big(\sigma^{-1}_s \nabla_{\nabla_\v X^{\ell^\eps}_{s-}}\sigma_s\big) (X^{\ell^\eps}_{s-}),\ \
 A_s:= \big(\sigma^{-1}_s \nabla_{\nabla_\v X^{\ell}_{s-}}\sigma_s\big) (X^{\ell}_{s-}).
\label{F3}\end{align}
Then
\begin{align*}
\eta^\eps_t:=\sum_{s\in[0,t]}\<\sigma^{-1}_s(X^{\ell^\eps}_{s-})\nabla_{\nabla_\v X^{\ell^\eps}_{s-}}\sigma_s(X^{\ell^\eps}_{s-})\Delta W^\beta_{\ell^\eps_{s}},\Delta W_{\ell^\eps_{s}}\>
=\sum_{s\in[0,t]}\<A^\eps_s \Delta W^\beta_{\ell^\eps_s}, \Delta W_{\ell^\eps_s}\>
\end{align*}
and
$$
\eta_t:=\sum_{s\in[0,t]}\<\sigma^{-1}_s(X^\ell_{s-})\nabla_{\nabla_\v X^\ell_{s-}}\sigma_s(X^\ell_{s-})\Delta W^\beta_{\ell_{s}},\Delta W_{\ell_{s}}\>
=\sum_{s\in[0,t]}\<A_s \Delta W^\beta_{\ell_s}, \Delta W_{\ell_s}\>.
$$
Letting $\|\cdot\|_p$ denote the $L^p$-norm with respsect to $\mathbb P,$ we obtain  for $p\in(1,2)$ that
$$
\|\eta^\eps_t-\eta_t\|_p\leq J_1(\vv)+J_2(\vv),
$$
where
\begin{align*}
J_1(\vv)&:=\sum_{s\in[0,t]}\|\<A_s \Delta W^\beta_{\ell^\eps_s}, \Delta W_{\ell^\eps_s}\>-\<A_s
\Delta W^\beta_{\ell_s}, \Delta W_{\ell_s}\>\|_p,\\
J_2(\vv)&:=\sum_{s\in[0,t]}\|\<(A^\eps_s-A_s)\Delta W^\beta_{\ell^\eps_s}, \Delta W_{\ell^\eps_s}\>\|_p.
\end{align*}
By H\"older's inequality, {\bf (A)} and Lemma \ref{Le2}, for $q=\frac{2p}{2-p}$ and some constants $C_1,C_2 $ depending on $p$ and $t$, we have
\begin{equation}\label{F11}\begin{split}
J_1(\vv)&\leq C_1 \sum_{s\in[0,t]}\|A_s\|_{2q}\|\Delta W_{\ell^\eps_s}\|_{2q}\|\Delta W^\beta_{\ell^\eps_s}-\Delta W^\beta_{\ell_s}\|_2\\
&\quad+C_1 \sum_{s\in[0,t]}\|A_s\|_{2q}\|\Delta W^\beta_{\ell_s}\|_{2q}\|\Delta W_{\ell^\eps_s}-\Delta W_{\ell_s}\|_2\\
&\leq C_2(1+|x|^m) \left(\sum_{s\in[0,t]}\|\Delta W_{\ell^\eps_s}\|^2_{2q}\right)^{1/2}
\left(\sum_{s\in[0,t]}\|\Delta W^\beta_{\ell^\eps_s}-\Delta W^\beta_{\ell_s}\|^2_2\right)^{1/2}\\
&\quad+C_2(1+|x|^m)\left(\sum_{s\in[0,t]}\|\Delta W^\beta_{\ell_s}\|^2_{2q}\right)^{1/2}
\left(\sum_{s\in[0,t]}\|\Delta W_{\ell^\eps_s}-\Delta W_{\ell_s}\|^2_2\right)^{1/2},
\end{split}\end{equation}
which converges to zero as $\vv\downarrow 0$ due to   the  argument   in the proof of Lemma \ref{Le1}, where $x$ is the initial point of $X_s^\ell$.

Similarly,   we   have
\begin{align}
J_2(\vv)\leq C\sum_{s\in[0,t]}\|A^\eps_s-A_s\|_{2}(\Delta\ell^\eps_s)^{1/2}(\Delta\lambda^\beta_{\ell^\eps_s})^{1/2}
\leq C\sup_{s\in[0,t]}\|A^\eps_s-A_s\|_{2}\label{F12}
\end{align}
which goes to zero as $\vv\downarrow 0$ by Lemma \ref{Le2}. By the dominated convergence theorem, and combining this with (\ref{F11}) and (\ref{F12}), we obtain
$$
\lim_{\vv\downarrow 0}I_3(\eps)= \frac{1}{\beta_{\ell_t}}\mE\left[f(X^{\ell}_t)\eta_t\right].
$$
Therefore,  the proof is finished by taking $\vv\downarrow 0$ in (\ref{F000}) and noting that due to (\ref{L1.4.3}),
$$\lim_{\vv\to 0} \nn_\v \mE f(X_t^{\ell_\vv}) =\nn_\v \mE f(X_t^\ell).$$
The proof is complete.
\end{proof}

\section{ Proof of Theorem 1.1}

Let $\mathcal B_b(\mR^d)$  be the set of all bounded measurable functions on $\mR^d$. Consider
$$
P_t^\ell f(x):=\mE f(X^{\ell}_t(x)),\ \ x\in\mR^d, t>0, f\in \mathcal B_b(\mR^d).
$$
We first prove the following gradient estimate on $P^\ell_t$.
\bt\label{Th00}
Assume $(\ref{Ass0})$ and $(\ref{Ass1})$.
For any $p>1$, there exists a constant $C >0$ such that
for any $f\in \mathcal B_b(\mR^d)$ and $x\in\mR^d$,  $t\in(0,1]$ with $\ell_{t}>0$,
\begin{align}
|\nabla P_t^\ell f(x)|\leq \ff {C(1+|x|^m)} {\ss{\ell_t\wedge 1}}(P^\ell_t|f|^p)^{1/p}(x).\label{Esw}
\end{align}
\et
\begin{proof} By   the Jensen inequality, we only need to prove  (\ref{Esw}) for   $p\in (1,2)$.
In this case, $q:= \ff p {p-1}>2.$

Simply denote $X_\cdot^\ell =X_\cdot^\ell (x)$. By (\ref{Ass1}), (\ref{L1.4.1}) and (\ref{L1.4.2}), there exists a constant $c>0$ such that for all $t\in(0,1]$,
\begin{align} \mE\left(\sup_{s\in [0,t]} |\nn_\v X_s^\ell|^{2q}\right)\leq c \e^{c\ell_t^q},\ \
\mE\left(\sup_{s\in [0,t]} |\si_s^{-1}( X_s^\ell)|^{2q}\right)\leq c (1+|x|^m)^{2q}\e^{c\ell_t^q}.\label{FG}
\end{align}
Define
$$
\tau:=\inf\{t: \ell_t\geq 1\}.
$$
Clearly,
\begin{align}
\ell_\tau\geq 1,\ \ \ell_{\tau-}\leq 1.\label{FG99}
\end{align}
Below we take
\begin{align}
\bb_t:=t\wedge \ell_{\tau},\ \  t\geq 0.\label{EQ}
\end{align}
By Theorem \ref{Main}, for any $\v\in\mR^d$ with $|\v|=1$, we have
\begin{align}
|\nabla_{\mathrm{v}}P_t^\ell f|\leq I_1+I_2+I_3,\label{FG0}
\end{align}
where
\begin{align*}
I_1&:=\frac{1}{\beta_{\ell_t}}\mE\left|f(X^{\ell}_t)
\int_0^t\sigma^{-1}_{s}(X^{\ell}_{s-})\left(\nabla_{\mathrm{v}} X^{\ell}_{s-}\right)\dif W^\beta_{\ell_{s}}\right|,\no\\
I_2&:=\frac{1}{\beta_{\ell_t}}\mE\left|f(X^{\ell}_t)
\int_0^t\mathrm{Tr}\big(\sigma^{-1}_s\nabla_{\nabla_\v X^{\ell}_{s-}}\sigma_s\big)(X^{\ell}_{s-}) \dif\beta_{\ell_s}\right|,\no\\
I_3&:=\frac{1}{\beta_{\ell_t}}\mE\left|f(X^{\ell}_t)\sum_{s\in[0,t]}\<A_s \Delta W^\beta_{\ell_s}, \Delta W_{\ell_s}\>\right|,
\end{align*}where $A_s$ is in (\ref{F3}).
For $I_1$, by Lemma \ref{Le0} and H\"older's  inequality, and using (\ref{FG}) and (\ref{LB}), we have
\begin{align*}
I_1&\leq \frac{1}{\beta_{\ell_t}}(P^\ell_t|f|^p)^{1/p}\left(\mE\left|
\int_0^t\sigma^{-1}_{s}(X^{\ell}_{s-})\left(\nabla_{\mathrm{v}} X^{\ell}_{s-}\right)\dif W^\beta_{\ell_{s}}\right|^q\right)^{1/q}\\
&\leq \frac{1}{\beta_{\ell_t}}(P^\ell_t|f|^p)^{1/p}\left(\mE\left(\int_0^t\left|\sigma^{-1}_{s}(X^{\ell}_{s-})
\left(\nabla_{\mathrm{v}} X^{\ell}_{s-}\right)\right|^2\dif \lambda^\beta_{\ell_{s}}\right)^{q/2}\right)^{1/q}\\
&\leq \frac 1{\beta_{\ell_t}}(P^\ell_t|f|^p)^{1/p}(\lambda^\beta_{\ell_t})^{\frac{1}{2}-\frac{1}{q}}\left(
\int_0^t\big(\mE \|\si_s^{-1}(X_{s-}^\ell)\|^{2q}\big)^{\ff 1 2}\big(\mE|\nabla_{\mathrm{v}} X^{\ell}_{s-}|^{2q}\big)^{\ff 1 2}\dif \lambda^\beta_{\ell_{s}}\right)^{\ff 1 q}\\
&\leq \frac{C_1(1+|x|^m)}{\beta_{\ell_t}}(P^\ell_t|f|^p)^{1/p}\bigg(\lambda^\beta_{\ell_t}\bigg)^{\ff 1 2-\ff 1 q}
\left(\int_0^t \e^{c\ell_{s-}^q}\dif\lambda^\beta_{\ell_s}\right)^{\ff 1 q}
\end{align*}
for some constant $C_1>0.$ By the choice (\ref{EQ}) of $\beta_t$ and the definition of $\lambda^\beta_{\ell_t}$, it is easy to see that
$$
\lambda^\beta_{\ell_t}=\int^{\ell_t}_0|\dot\beta_s|^2\dif s=\int^{\ell_t}_01_{s\leq \ell_{\tau}}\dif s=\ell_t\wedge\ell_{\tau}=\beta_{\ell_t},
$$
and  by (\ref{FG99}),
$$
\int_0^t \e^{c\ell_{s-}^q}\dif\lambda^\beta_{\ell_s}=\int_0^t \e^{c\ell_{s-}^q}\dif\beta_{\ell_s}=\sum_{s\in[0,t]} \e^{c\ell_{s-}^q}\Delta\beta_{\ell_s}\leq \e^c\beta_{\ell_t}.
$$
Therefore,
\begin{align}
I_1\leq\frac{C_1(1+|x|^m) \e^{c/q}}{\sqrt{\beta_{\ell_t}}}(P^\ell_t|f|^p)^{1/p}
\stackrel{(\ref{FG99})}{\leq}\frac{C_1(1+|x|^m) \e^{c/q}}{\sqrt{\ell_t\wedge 1}}(P^\ell_t|f|^p)^{1/p}.\label{FG2}
\end{align}
Similarly, there exists a constant $C_2>0$ such that
\begin{equation}\label{FG3}\begin{split}
I_2&\leq\frac{1}{\beta_{\ell_t}}(P^\ell_t|f|^p)^{1/p}\left(\mE\left|
\int_0^t\mathrm{Tr}\big(\sigma^{-1}_s\nabla_{\nabla_\v X^{\ell}_{s-}}\sigma_s\big)(X^{\ell}_{s-})\dif \beta_{\ell_{s}}\right|^q\right)^{1/q}\\
&\leq\frac{1}{\beta_{\ell_t}}(P^\ell_t|f|^p)^{1/p}(\beta_{\ell_t})^{1-\frac{1}{q}}\left(
\int_0^t\mE|\mathrm{Tr}\big(\sigma^{-1}_s\nabla_{\nabla_\v X^{\ell}_{s-}}\sigma_s\big)(X^{\ell}_{s-})|^q\dif \beta_{\ell_{s}}\right)^{1/q}\\
&\leq C_2(1+|x|^m)(P^\ell_t|f|^p)^{1/p}(\beta_{\ell_t})^{-\frac{1}{q}}\left(\int_0^t\e^{c\ell_{s-}^q}\dif \beta_{\ell_{s}}\right)^{1/q}\\
&\leq C_2(1+|x|^m)\e^{c/q}(P^\ell_t|f|^p)^{1/p}.\end{split}\end{equation}
 Finally,  noting that $\DD\bb_{\ell_s}= \DD\ell_s 1_{[0,\tau]}(s)$, we have, for some constant  $C_3>0$,
\begin{align*}
I_3&\leq\frac{1}{\beta_{\ell_t}}(P^\ell_t|f|^p)^{1/p}\left\|\sum_{s\in[0,t]}\<A_s \Delta W^\beta_{\ell_s}, \Delta W_{\ell_s}\>\right\|_q\\
&\leq\frac{1}{\beta_{\ell_t}}(P^\ell_t|f|^p)^{1/p}\sum_{s\in[0,t]}\|\<A_s \Delta W^\beta_{\ell_s}, \Delta W_{\ell_s}\>\|_q\\
&\leq\frac{1}{\beta_{\ell_t}}(P^\ell_t|f|^p)^{1/p}\sum_{s\in[0,t]}\|A_s\|_{2q}\|\Delta W^\beta_{\ell_s}\|_{4q}\|\Delta W_{\ell_s}\|_{4q}\\
&\leq\frac{C_3(1+|x|^m)}{\beta_{\ell_t}}(P^\ell_t|f|^p)^{1/p}\sum_{s\in[0,t]}|\Delta\ell_s|^{1/2}|\Delta\lambda^\beta_{\ell_s}|^{1/2}\e^{c\ell_{s-}^q/(2q)}\\
&\leq\frac{C_3(1+|x|^m)}{\beta_{\ell_t}}(P^\ell_t|f|^p)^{1/p}\sum_{s\in[0,t]}|\Delta\ell_s|^{\frac{1}{2}}|\Delta\beta_{\ell_s}|^{\frac{1}{2}}\e^{c/(2q)}\\
&= \ff{C_3\e^{c/(2q)} (1+|x|^m)}{\bb_{\ell_t}}(P^\ell_t|f|^p)^{1/p} \sum_{s\in [0,t\land\tau]} \DD\ell_s  \\
&=C_3(1+|x|^m) \e^{c/(2q)}(P^\ell_t|f|^p)^{1/p}.\end{align*}
Combining this with (\ref{FG2}), (\ref{FG3}) and (\ref{FG0}), we obtain (\ref{Esw}).
\end{proof}

Now we are in a position to give

\vspace{2mm}

\begin{proof}[Proof of Theorem 1.1]
By taking regular conditional probability with respect to $S$, we have
\begin{align}
P_tf(x):=\mE f(X_t(x))=\mE\big(\mE f(X^{\ell}_t(x))|_{\ell=S}\big)=\mE\big(P^S_tf(x)\big).\label{F33}
\end{align}

(i) For $R>0$ and $\ell\in\mS$, define
$$
\tau(\ell):=\inf\{t: \ell_t\geq R\}.
$$
If we choose $\beta_t=t\wedge\tau(\ell)$ in (\ref{EE1}), then by (\ref{F33}) and (\ref{EE1}), we have, for any $\v\in\mR^d$ with $|\v|=1$,
\begin{align*}
\nabla_{\mathrm{v}} P_tf&=I_1+I_2+I_3,
\end{align*}
where
\begin{align*}
I_1&:=\mE\left(\mE\Bigg[f(X^\ell_t)\frac{1}{\beta_{\ell_t}}\int^t_0
\<\sigma^{-1}_s(X^\ell_{s-})\nabla_\v X^\ell_{s-},\dif W^\beta_{\ell_s}\>\Bigg]_{\ell=S}\right),\\
I_2&:=\mE\left(\mE\Bigg[f(X^\ell_t)\frac{1}{\beta_{\ell_t}}\int^t_0\mathrm{Tr}\big(\sigma^{-1}_s\nabla_{\nabla_\v
X^\ell_{s-}}\sigma_s\big) (X^\ell_{s-}) \dif \beta_{\ell_s}\Bigg]_{\ell=S}\right),\no\\
I_3&:=\mE\left(\mE\Bigg[f(X^\ell_t)\frac{1}{\beta_{\ell_t}}
\sum_{s\in[0, t]}\<\sigma^{-1}_s(X^\ell_{s-})\nabla_{\nabla_\v X^\ell_{s-}}\sigma_s(X^\ell_{s-})\Delta W^\beta_{\ell_{s}},\Delta W_{\ell_{s}}\>\Bigg]_{\ell=S}\right).
\end{align*}
As shown in the proof of  Theorem \ref{Th00}, it is clear that when $\mE S_{t\land\tau}^{-1/2}<\infty,$ $I_1, I_2$ and  $I_3$ are well defined. Noticing that
$$
W^\beta_{\ell_s}=W_{\ell_s\wedge \ell_{\tau(\ell)}}=W_{\ell_{s\wedge \tau(\ell)}}
$$
and
$$
\Delta W^\beta_{\ell_s}=\Delta W_{\ell_{s\wedge \tau(\ell)}}=W_{\ell_s}1_{[0,\tau(\ell)]}(s),
$$
we have
\begin{align*}
I_1&=\mE\Bigg[f(X_t)\frac{1}{S_{t\wedge\tau}}\int^{t\wedge\tau}_0
\<\sigma^{-1}_s(X_{s-})\nabla_\v X_{s-},\dif W_{S_s}\>\Bigg],\\
I_2&=\mE\Bigg[f(X_t)\frac{1}{S_{t\wedge\tau}}\int^{t\wedge\tau}_0\mathrm{Tr}\big(\sigma^{-1}_s\nabla_{\nabla_\v X_{s-}}\sigma_s\big) (X_{s-}) \dif S_s\Bigg],
\end{align*}
and
\begin{align*}
I_3&=\mE\Bigg[f(X_t)\frac{1}{S_{t\wedge\tau}}\sum_{s\in[0, t\wedge\tau]}\<\sigma^{-1}_s(X_{s-})\nabla_{\nabla_\v X_{s-}}
\sigma_s(X_{s-})\Delta W_{S_{s}},\Delta W_{S_{s}}\>\Bigg)\Bigg]\\
&=\mE\Bigg[f(X_t)\frac{1}{S_{t\wedge\tau}}\int^{t\wedge\tau}_0\!\!\!\int_{\mR^d}\<\sigma^{-1}_s(X_{s-})\nabla_{\nabla_\v X_{s-}}\sigma_s(X_{s-})x,x\>N(\dif s,\dif x)\Bigg],
\end{align*}
where the random measure $N$ is defined by (\ref{FF99}). Thus, the formula (\ref{EE11}) is proven.

(ii) By the Markov property, it suffices to prove  (\ref{FG88}) for $t\in(0,1]$. Then the estimate (\ref{FG88}) follows from Theorem 4.1 and (\ref{F33}) by using the H\"older inequality. So, it suffices to prove (\ref{Ess}) for $t\in (0,1]$ in the case that $B(u)\ge c u^{\alpha/2}, u\ge u_0$, where
$\alpha\in(0,2)$ and $u_0\geq 0$. In this case, for any $\gamma>0$ we may find a constant $C>0$ such that for all $t>0$ (see e.g. \cite[page 298]{GRW} for the formula of $\mE S_t^{-\gamma}$),
\begin{align*} \mE\ff 1{ S_t^\gamma\land 1 }&\le 1 +\mE \ff 1 {S_t^\gamma} =1+ \ff 1{\GG(\gamma)} \int_0^\infty u^{\gamma-1}\e^{-tB(u)}\dif u\\
&\le 1+ \ff 1{\GG(\gamma)}\bigg( \int_{u_0}^\infty u^{\gamma-1}\e^{-ctu^{\aa/2}}\dif u + \int_0^{u_0} u^{\gamma-1}\d u\bigg)\\
&\le 1+ \ff {u_0^\gamma}{\gamma\GG(\gamma)}  +\ff1 {\GG(\gamma)}\int_0^\infty u^{\gamma-1}\e^{-ctu^{\aa/2}}\dif u\le \ff C{1\land t^{2\gamma/\aa}}.
\end{align*}
  Therefore, (\ref{Ess}) follows from (\ref{FG88}).
\end{proof}

{\bf Acknowledgements.} We would like to   thank Professor Zhen-Qing Chen for    useful conversations as well as the referees for helpful comments.

\end{document}